\documentclass[11 pt,a4paper,leqno]{article}
\usepackage{latexsym,amsmath,amscd,amsfonts}

\RequirePackage{amsmath}
\RequirePackage{amssymb}
\RequirePackage{ifthen}
\RequirePackage{theorem}
\RequirePackage{pslatex}

\title{\sc The $q$-analogue of the wild fundamental group (II)}

\author{J.-P. Ramis\footnote{Institut de France (Acad\'emie des Sciences)},\;
J. Sauloy
\footnote{Laboratoire Emile Picard, CNRS UMR 5580, U.F.R. M.I.G.,
118, route de Narbonne, 31062 Toulouse CEDEX 4}}

\date{}

\def\sq{\sigma_q}

\def\C{{\mathbf C}}

\def\Z{{\mathbf Z}}

\def\N{{\mathbf N}}

\def\ii{{\text{i}}}
\def\F{{\mathcal{F}}}

\def\M{{\mathcal{M}}}
\def\Ma{{\text{Mat}}}
\def\Co{{\mathcal{C}}}
\def\D{{\mathcal{D}_{q,K}}}
\def\gr{{\text{gr}}}

\def\Eq{{\mathbf{E}_{q}}}
\def\EE{{\mathcal{E}}}
\def\G{{\mathfrak{G}_{A_{0}}}}
\def\g{{\mathfrak{g}_{A_{0}}}}

\def\Ra{{\C\{z\}}}
\def\Ka{{\C(\{z\})}}

\def\Raq(d){{\C\{\xi\}_{q,(\delta)}}}
\def\Kaq(d){{\C(\{\xi\})_{q,(\delta)}}}

\def\Rf{{\C[[z]]}}
\def\Kf{{\C((z))}}
\def\Rw{{\mathcal{O}(\C^{*})}}
\def\Kw{{\mathcal{M}(\C^{*})}}
\def\Rwg{{\mathcal{O}(\C^{*},0)}}
\def\Kwg{{\mathcal{M}(\C^{*},0)}}

\def\Der{{\dot{\Delta}}}
\def\Derc{{\Der_{\overline{c}}}}
\def\Derdc{{\Der^{(\delta)}_{\overline{c}}}}

\def\Lsca{{LS_{\overline{c},a}}}
\def\Lsda{{LS_{\overline{d},a}}}
\def\homa{{\hat{\omega}^{(0)}_{a}}}

\def\diag{\text{diag}}

\def\Sp{\text{Sp}}

\def\St{\mathfrak{St}}
\def\st{\mathfrak{st}}

\newtheorem{thm}{Theorem}[section]
\newtheorem{lemma}[thm]{Lemma}
\newtheorem{prop}[thm]{Proposition}
\newtheorem{cor}[thm]{Corollary}

\def \Pr {\textsl{Proof. - }}
\def\Ex{\noindent\textbf{Example.~}}

\def\Rem{\noindent\textbf{Remark.~}}
\def\Rems{\noindent\textbf{Remarks. \\}}

\begin{document}

\maketitle

\bigskip \hrule \bigskip

\begin{abstract}
{\small 
In \cite{RS1}, we defined $q$-analogues of alien derivations
and stated their basic properties. In this paper, we prove 
the density theorem and the freeness theorem announced in 
\emph{loc. cit.}.}
\end{abstract}

\bigskip \hrule 

\tableofcontents

\bigskip \hrule \bigskip



\section{Introduction}


\subsection{The problem}

In this paper we shall continue the study of the local meromorphic 
classification of $q$-difference modules. In \cite{RSZ} we gave such 
a classification in Birkhoff style, using normal forms and index theorems; 
this classification is complete in the ``integral slope case". (One could 
extend it to the general case using some results of \cite{vdPR}.) \\
 
In \cite{RS1} we introduced a new approach of the classification, using
a ``fundamental group'' and its finite dimensional representations, 
in the style of the Riemann-Hilbert correspondence for linear differential 
equations. At some abstract level, such a classification is well known: 
the fundamental group is the tannakian Galois group of the tannakian 
category of local meromorphic $q$-modules. But we wanted more information: 
our essential aim was to get a \emph{smaller} fundamental group which is 
Zariski dense in the tannakian Galois group and to describe it 
{\it explicitly}, in the spirit of what was done by the first author 
for the differential case \cite{RamInv3}. \\

In \cite{RS1} we built a family of elements of the Lie algebra of 
the tannakian group, the $q$-alien derivations, we achieved our program 
for the one-level case and we announced the main results in the general case.
The aim of the present paper is to give some proofs omitted in \cite{RS1} 
for the multi-level case. We will finally give a more precise algebraic 
formulation of our results in \cite{RS3}, which will end the series.


\subsection{Contents of the paper}

General notations and conventions are explained in the next paragraph
\ref{subsection:generalnotations}. In section \ref{section:analyticQED}, 
we recall basic properties of the category $\EE^{(0)}_{1}$ of linear 
analytic $q$-difference equations with integral slopes, and the structure 
and action of its Galois group $G^{(0)}_{1}$. In section 
\ref{section:wildgroup}, we recall the unipotent structure of the Stokes 
subgroup $\St$ of $G^{(0)}_{1}$, and the construction (taken from 
\cite{RS1}) of elements of the Lie algebra $\st$ of $\St$, the 
so-called \emph{$q$-alien derivations}. Our ``$q$-analogue of the
wild fundamental group'' is the Lie subalgebra they generate. We then
prove in \ref{subsection:density} and \ref {subsection:freeness} our
main results: density and a freeness property of the $q$-alien derivations. 
In section \ref{section:conclusion}, we summarize what remains to 
be solved, and will be the contents of \cite{RS3}. \\

The paper is written so as to be read widely independently
of \cite{RS1} - granted the reader is willing to take on faith 
some key points. The principle of the proofs is almost purely 
tannakian, but we have stated explicitly the underlying methods 
and prerequisites. Moreover, they are described in a concrete, 
computational form (with a systematic use of matrices). Since 
neither $q$-difference equations, nor even tannakian methods are 
so popular, this may help the reader to get familiarized with
either domain. Note that, since we heavily rely on transcendental 
tools, the methods here are, to a large extent, independent of those 
of M. van der Put and his coauthors.


\subsection{General notations}
\label{subsection:generalnotations}

The notations are the same as in \cite{RS1}. Here are the
most useful ones. \\

We let $q \in \C$ be a complex number with modulus $|q| > 1$.
We write $\sq$ the $q$-dilatation operator, so that, for any map $f$ 
on an adequate domain in $\C$, one has: $\sq f(z) = f(qz)$.
Thus, $\sq$ defines a ring automorphism in each of the following rings:
$\Ra$ (convergent power series), $\Rf$ (formal power series), 
$\Rw$ (holomorphic functions over $\C^{*}$), $\Rwg$ (germs at $0$ of 
elements of $\Rw$). Likewise, $\sq$ defines a field automorphism in each 
of their fields of fractions: $\Ka$ (convergent Laurent series), $\Kf$ 
(formal Laurent series over), $\Kw$ (meromorphic functions over $\C^{*}$), 
$\Kwg$ (germs at $0$ of elements of $\Kw$) \\

The $\sq$-invariants elements of $\Kw$ can be considered as meromorphic
functions on the quotient Riemann surface $\Eq = \C^{*}/q^{\Z}$. Through 
the mapping $x \mapsto z = e^{2 \ii \pi x}$, the latter is identified to
the complex torus $\C/(\Z + \Z \tau)$, where $q = e^{2 \ii \pi \tau}$.
Accordingly, we shall identify the fields $\Kw^{\sq}$ and $\M(\Eq)$. 
We shall write $a \mapsto \overline{a}$ the canonical projection map
$\pi: \C^{*} \rightarrow \Eq$ and 
$[c;q] = \pi^{-1}\left(\overline{c}\right) = c q^{\Z}$
(a discrete logarithmic $q$-spiral). \\

Last, we shall have use for the function $\theta \in \Rw$, a Jacobi Theta 
function such that $\sq \theta = z \theta$ and $\theta$ has simple zeroes 
along $[-1;q]$. One then puts $\theta_{c}(z) = \theta(z/c)$, so that
$\theta_{c} \in \Rw$ satisfies $\sq \theta_{c} = (z/c) \theta_{c}$ and 
$\theta_{c}$ has simple zeroes along $[-c;q]$.


 
\section{Linear analytic $q$-difference equations}
\label{section:analyticQED}

A linear analytic (resp. formal) $q$-difference equation (implicitly:
at $0 \in \C$) is an equation:
\begin{equation}
\label{equation:qED}
\sq X = A X, 
\end{equation}
where $A \in GL_{n}(\Ka)$ (resp. $A \in GL_{n}(\Kf)$). There is 
an intrinsic description as a ``$q$-difference module $M_{A}$'',
which runs as follows. We consider the operator $\Phi_{A}$ on
$\Ka^{n}$ which maps a column vector $X$ to $A^{-1} \sq X$.
This can be abstracted as a finite dimensional $\Ka$-vector
space $V$ endowed with a so-called ``$\sq$-linear automorphism'' $\Phi$.
A $q$-difference module is such a pair $M = (V,\Phi)$. Here, we have
$M_{A} = (\Ka^{n},\Phi_{A})$. \\

We shall here stick to the matrix model and, for all practical
purposes, the reader may identify the equation (\ref{equation:qED}),
the matrix $A$ and the $q$-difference module $M_{A}$ with each other.
For instance, we call \emph{solution of $A$, or of (\ref{equation:qED}), 
or of $M_{A}$} in some extension $K$ of $\Ka$, on which $\sq$ operates, 
a column vector $X \in K^{n}$ such that $\sq X = A X$. The \emph{underlying 
space} of $A \in GL_{n}(\Ka)$ is $\Ka^{n}$.


\subsection{Description of the tannakian structure}
\label{subsection:tannakianstructure}

We now proceed to describe the \emph{tannakian category of analytic 
$q$-difference equations} $\EE^{(0)}$. There is a similar description 
for the corresponding formal category $\hat{\EE}^{(0)}$. The objects 
of  $\EE^{(0)}$ are linear analytic $q$-difference equations 
(\ref{equation:qED}). A morphism from $A \in GL_{n}(\Ka)$
to $B \in GL_{p}(\Ka)$ is a matrix $F \in M_{p,n}(\Ka)$ such that:
\begin{equation}
\label{equation:morphisme}
(\sq F) A = B F. 
\end{equation}
This just means that $F$ sends any solution $X$ of $A$ to a solution
$F X$ of $B$. One can check that $\EE^{(0)}$ is an abelian category.
Indeed, it is the category of finite length left modules over the euclidean
non-commutative ring $\D$ of $q$-difference operators over $K = \Ka$. \\

The abelian category $\EE^{(0)}$ is endowed with a tensor structure. 
The tensor product $A_{1} \otimes A_{2}$ of two objects
(resp. the tensor product $F_{1} \otimes F_{2}$ of two morphisms)
is just the Kronecker product of the matrices; of course, we must
define a consistent way of identifying $\C^{n} \otimes \C^{p}$
with $\C^{np}$, or $\Ka^{n} \otimes \Ka^{p}$ with $\Ka^{np}$ (see, 
for instance \cite{JSGAL}). \\

The unit object $\underline{1}$ (which is neutral for the tensor product) 
is the matrix $(1) \in GL_{1}(\Ka) = \Ka^{*}$, with underlying space $\Ka$.
The object $\underline{1}$ of course corresponds to the ``trivial''
equation 
\footnote{In differential Galois theory, the matrix $A$ of a system
is in $M_{n}(\Ka)$ (rather than in $GL_{n}$), the trivial equation is 
$x' = 0$, etc. The theory of $q$-difference equations rather has a 
multiplicative character}
$\sigma x = x$.
One easily checks that the space $Hom(\underline{1},A)$ of morphisms from 
$\underline{1}$ to $A$ is exactly the space of solutions of $A$ in $\Ka$,
or, equivalently, the space of fixed points of $\Phi_{A}$ in $\Ka^{n}$.
We shall write $\Gamma(A)$ or $\Gamma(M_{A})$ that space, as it is similar 
to a space of global sections (and, indeed, can be realised as such, see
\cite{JSFIB}). \\

The characterization (\ref{equation:morphisme}) of morphisms can itself
be seen as a $q$-difference equation $\sq F = B F A^{-1}$. This means
that there is an ``internal Hom'' object, which can be described in
the following way. Consider the linear operator $F \mapsto B F A^{-1}$ 
on the vector space $M_{p,n}(\Ka)$. Through identification of $M_{p,n}(\Ka)$ 
with $\Ka^{np}$, this operator is described by a matrix in $GL_{np}(\Ka)$,
which yields the desired equation. We shall write $\underline{Hom}(A,B)$
the corresponding object. Thus, one gets:
\begin{equation}
\label{equation:hominterne1}
\Gamma(\underline{Hom}(A,B)) \simeq Hom(\underline{1},\underline{Hom}(A,B)) 
\simeq Hom(A,B).
\end{equation}
Actually, this is a special case of the following canonical isomorphism::
\begin{equation}
\label{equation:hominterne2}
Hom(A,\underline{Hom}(B,C)) \simeq Hom(A \otimes B,C).
\end{equation}
The reader will check that the object $\underline{Hom}(A,\underline{1})$ 
has the following description. The underlying space is $M_{1,n}(\Ka)$, 
which we identify with $\Ka^{n}$. The corresponding matrix for the linear 
operator $F \mapsto F A^{-1}$ is the \emph{contragredient matrix} 
$A^{\vee} = {}^{t}A^{-1}$. We call the object $A^{\vee}$ the \emph{dual} 
of the object $A$. From this, we get yet another construction of the 
internal Hom:
\begin{equation}
\label{equation:hominterne3}
\underline{Hom}(A,B) \simeq A^{\vee} \otimes B.
\end{equation}
We summarize these properties by saying that $\EE^{(0)}$ is
a tannakian category. This is halfway to showing that it is
(isomorphic to) the category of representations
of a proalgebraic group, our hoped for Galois group. To get
further, one needs a \emph{fiber functor on $\EE^{(0)}$}. 
This was defined and, to some extent, studied in full generality 
in \cite{JSFIL}, \cite{JSSTO} and \cite{RS1}. However, for our
strongest results, we need to restrict to the case of \emph{integral 
slopes}.


\subsection{Equations with integral slopes}
\label{subsection:integralslopes}

In \cite{JSFIL}, one defined the Newton polygon of a $q$-difference
module (analytic or formal). This consists in slopes
\footnote{Note that here, as in \cite{RS1}, we have changed
the definitions of slopes. Those used here are the opposites
of those used in \cite{JSFIL}, \cite{RSZ1}and \cite{JSSTO}.}
$\mu_{1} < \cdots < \mu_{k}$ (rational numbers) together with ranks 
(or multipicities) $r_{1},\ldots,r_{k}$ (positive integers). 
We shall say that a module is \emph{pure isoclinic} if it has 
only one slope and that it is \emph{pure}
\footnote{Here again, starting with \cite{RS1}, we changed our
terminology: we now call pure isoclinic (resp. pure) what was 
previously called pure (resp. tamely irregular). The latter are 
called \emph{split modules} in \cite{vdPR}.}
if it is a direct sum of pure isoclinic modules. We call
\emph{fuchsian} a pure isoclinic module with slope $0$.
The Galois theory of fuchsian modules was studied in \cite{JSGAL}. 
Pure modules are irregular objects without \emph{wild monodromy}, 
as follows from \cite{RSZ}, \cite{JSSTO} and \cite{RS1}. \\
 
The tannakian subcategory of $\EE^{(0)}$ made up of pure modules
is called $\EE^{(0)}_{p}$. Modules with integral slopes also form
tannakian subcategories, which we write $\EE^{(0)}_{1}$ and 
$\EE^{(0)}_{p,1}$. \textbf{\emph{From now on, we restrict to the case
of integral slopes}}. Our category of interest is therefore 
$\EE^{(0)}_{1}$ and we shall now start its description. \\

Any equation in $\EE^{(0)}_{1}$ can be written in the following
\emph{standard form}:
\begin{equation} 
\label{equation:formestandard}
A = 
\begin{pmatrix}
z^{\mu_{1}} A_{1}  & \ldots & \ldots & \ldots & \ldots \\
\ldots & \ldots & \ldots  & U_{i,j} & \ldots \\
0      & \ldots & \ldots   & \ldots & \ldots \\
\ldots & 0 & \ldots  & \ldots & \ldots \\
0      & \ldots & 0       & \ldots & z^{\mu_{k}} A_{k}    
\end{pmatrix},
\end{equation}
where the slopes $\mu_{1} < \cdots < \mu_{k}$ are integers,
$r_{i} \in \N^{*}$, $A_{i} \in GL_{r_{i}}(\C)$ ($i = 1,\ldots,k$) 
(those $\mu_{i}$ and $r_{i}$ make up the Newton polygon of $A$) and:
$$
\forall i,j \text{~s.t.~} 1 \leq i < j \leq k \;,\;
U_{i,j} \in  \Ma_{r_{i},r_{j}}(\Ka).
$$
We actually can, and will, require the blocks $U_{i,j}$ to have 
all their coefficients in $\C[z,z^{-1}]$. Then any morphism
$F: A \rightarrow B$ between two matrices in standard form
is easily seen to be meromorphic at $0$ (by definition) and
holomorphic all over $\C^{*}$; this is because the equation
$\sq F = B F A^{-1}$ allows one to propagate the regularity
near $0$ to increasing neighborhoods. \\

We moreover say that $A$ is in \emph{polynomial} standard form 
if each block $U_{i,j}$ with $1 \leq i < j \leq k$ has coefficients 
in $\sum\limits_{\mu_{i} \leq d < \mu_{j}} \C z^{d}$. It was
proved in \cite{RSZ} that any object in $\EE^{(0)}_{1}$ is
analytically equivalent to one written in polynomial standard 
form (in essence, this is due to Birkhoff and Guenther). Last, 
we say that $A$ is in \emph{normalized} standard form is  if all 
the eigenvalues of all the blocks $A_{i}$ are in the fundamental 
annulus $\{z \in \C^{*} \mid 1 \leq |z| < |q|\}$. Any standard form 
can be normalized through shearing transformations. Note that
polynomial standard form is stable under tensor product, while 
normalized standard form is not. \\

The standard form (\ref{equation:formestandard}) above expresses 
the existence of a \emph{filtration by the slopes} (\cite{JSFIL}). 
The functoriality of the filtration moreover entails that a morphism
$F:A \rightarrow B$ is also upper triangular (by blocks) in the 
following sense: if the slopes of 
$B \in GL_{p}(\Ka)$ are $\nu_{1} < \cdots < \nu_{l}$, 
with ranks $s_{1} < \cdots < s_{l}$, then the morphism $F \in M_{p,n}(\Ka)$ 
from $A$ to $B$ has only non null blocks $F_{i,j} \in M_{s_{j},r_{i}}(\Ka)$, 
$1 \leq i \leq k$, $1 \leq j \leq l$ for $\nu_{j} \leq \mu_{i}$. \\

To the matrix $A$ and module $M = M_{A}$ is associated the graded
module $\gr M = M_{0} = M_{A_{0}}$ with block diagonal matrix:
\begin{equation} 
\label{equation:formestandardpure}
A_{0} = 
\begin{pmatrix}
z^{\mu_{1}} A_{1}  & \ldots & \ldots & \ldots & \ldots \\
\ldots & \ldots & \ldots  & 0 & \ldots \\
0      & \ldots & \ldots   & \ldots & \ldots \\
\ldots & 0 & \ldots  & \ldots & \ldots \\
0      & \ldots & 0       & \ldots & z^{\mu_{k}} A_{k}    
\end{pmatrix},
\end{equation}
The graded module $M_{0}$ is the direct sum 
$P_{1} \oplus \cdots \oplus P_{k}$, where each module $P_{i}$ 
is pure of rank $r_{i}$ and slope $\mu_{i}$ and corresponds
to the matrix $z^{\mu_{i}} A_{i}$. 
The functor $M \leadsto \gr M$ also acts on morphisms. To $F$,
it associates $F_{0}$ which has the same diagonal blocks as $F$,
that is, $(F_{0})_{i,j} = F_{i,j}$ if $\mu_{i} = \nu_{j}$. But
all the $(F_{0})_{i,j}$ such that $\mu_{i} \neq \nu_{j}$ are null. \\

By formalisation, \emph{i.e.} base change $\Ka \to \Kf$, the slope
filtration splits and the functor $\gr$ becomes isomorphic to the
identity functor. In matrix terms, this translates as follows.
There is a unique isomorphism $F: A_{0} \rightarrow A$ with formal
components $F_{i,j} \in M_{r_{i},r_{j}}(\Kf)$ (for $1 \leq i,j \leq k$)
and the following shape:
\begin{equation} 
\label{equation:automorphisme}
F = 
\begin{pmatrix}
I_{r_{1}} & \ldots & \ldots & \ldots & \ldots \\
\ldots & \ldots & \ldots  & F_{i,j} & \ldots \\
0      & \ldots & \ldots   & \ldots & \ldots \\
\ldots & 0 & \ldots  & \ldots & \ldots \\
0      & \ldots & 0       & \ldots & I_{r_{k}}     
\end{pmatrix}.
\end{equation}
To express that a matrix has such a shape and coefficients 
in some domain $K$, we shall write $F \in \G(K)$. Thus, $\G$
is a unipotent algebraic subgroup of the linear group and it
can be realised above any field $K$: in the above case, one
has $F \in \G(\Kf)$. For further use, we also give a notation
for the corresponding Lie algebra $\g$. An element $f \in \g(K)$
has the shape:
\begin{equation} 
\label{equation:endomorphisme}
\begin{pmatrix}
0_{r_{1}} & \ldots & \ldots & \ldots & \ldots \\
\ldots & \ldots & \ldots  & f_{i,j} & \ldots \\
0      & \ldots & \ldots   & \ldots & \ldots \\
\ldots & 0 & \ldots  & \ldots & \ldots \\
0      & \ldots & 0       & \ldots & 0_{r_{k}}     
\end{pmatrix},
\end{equation}
where $0_{r}$ is the null $r \times r$ matrix and where each
$f_{i,j} \in M_{r_{i},r_{j}}(K)$. \\

We shall denote $\hat{F}_{A}$ the unique $F$ mentioned above.
Its blocks can be characterised as the unique formal solutions
to the following recursive equations:
\begin{equation}
\label{equation:calculrecursifhat(F)}
\forall 1 \leq i < j \leq k \;,\;
\sq F_{i,j} z^{\mu_{j}} A_{j} - z^{\mu_{i}} A_{i} F_{i,j} =
\sum_{i < l < j} U_{i,l} F_{l,j} + U_{i,j}.
\end{equation}
There are usually no analytic solutions (that is, with coefficients
in $\Ka$) for equations (\ref{equation:calculrecursifhat(F)}). (The
existence of analytic solutions is equivalent to $M_{A}$ being pure.)
There are, however, meromorphic solutions, to be considered as
\emph{resummations} of the formal solution $\hat{F}_{A}$ (section 
\ref{subsection:Stokesaction}). \\

The graded counterpart $F_{0}$ of $F = \hat{F}_{A}$ satisfies simpler 
equations. From the above description, we know that $(F_{0})_{i,j} = 0$ 
for any $i,j$ such that $\mu_{i} \neq \mu_{j}$, that is, if $i \neq j$; 
if $i = j$:
$$
\sq (F_{0})_{i,i} z^{\mu_{i}} A_{i} = z^{\mu_{i}} A_{i} (F_{0})_{i,i}.
$$
This implies that $\sq (F_{0})_{i,i} A_{i} = A_{i} (F_{0})_{i,i}$,
and it then follows from \cite{JSGAL} that the coefficients of $F_{0}$
are Laurent polynomials (elements of $\C[z,z^{-1}]$); if moreover $A$
is in normalized standard form, then these coefficients are in $\C$.


\subsection{Description of the fiber functor}
\label{subsection:fiberfunctor}

In Tannaka theory, the Galois group is defined as the group
of tensor automorphisms of a fiber functor. We now describe
a fiber functor on $\EE^{(0)}_{1}$. There is actually a whole
family of these, indexed by $\C^{*}$, and one can therefore 
define a Galois groupo\"id (\cite{RS1}). Here, we shall first
choose an arbitrary basepoint $a \in \C^{*}$. As a consequence,
some constructions of \ref{subsection:Stokesaction} 
will be valid for most equations, but not all. This means that,
to study a particular equation, one has to choose a basepoint
compatible with it, which will be seen to be a generically
true condition. \\

The fiber functor $\homa$ goes from $\EE^{(0)}_{1}$ to the category 
of finite dimensional $\C$-vector spaces. On the side of objects, 
to each matrix $A \in GL_{n}(\Ka)$ and module $M_{A}$, it associates
the space $\homa(A) = \C^{n}$. On the side of
morphisms, to $F: A \rightarrow B \in GL_{p}(\Ka)$, it associates
$F_{0}(a): \C^{n} \rightarrow \C^{p}$. (The dimensions are right
and it follows from the last remark in \ref{subsection:integralslopes}
that $F_{0}(a)$ is well defined). \\

Apart from functoriality, the properties of $\homa$ 
which make it a fiber functor are the following: it is exact, faithful
and $\otimes$-compatible. The latter means that, for any $A,B$, 
the natural map 
$t_{A,B}: \homa(A) \otimes \homa(B)
\rightarrow \homa(A \otimes B)$ 
is an isomorphism. \\

We now define the Galois group of $\EE^{(0)}_{1}$ (at base point $a$)
as $G^{(0)}_{1} = Aut^{\otimes}(\homa)$. It would
be more rigorous to write explicitly the index $a$ indicating 
the basepoint, but this would make the notation heavier without 
true necessity. An element of the group 
$Aut^{\otimes}(\homa)$ is, 
by definition, a natural transformation
$g: A \leadsto g(A) \in 
GL\left(\homa(A)\right) = GL_{n}(\C)$, subject
to the following conditions:
\begin{enumerate}
\item{Functoriality: for any morphism $F: A \rightarrow B$, one has
$g(B) \circ F_{0}(a) = F_{0}(a) \circ g(A)$. Thus, the 
following diagram is commutative:
\begin{equation*}
\begin{CD}
\homa(A) @>F_{0}(a)>> \homa(B) \\
@V{g(A)}VV                          @VV{g(B)}V      \\
\homa(A) @>F_{0}(a)>> \homa(B)
\end{CD}
\end{equation*}
}
\item{Tensor compatibility: for any objects $A,B$, up to the natural
identifications, one has an equality
$g(A \otimes B) = g(A) \otimes g(B)$. Thus, the 
following diagram is commutative:
\begin{equation*}
\begin{CD}
\homa(A) \otimes \homa(B)
@>t_{A,B}>>
\homa(A \otimes B) \\
@V{g(A) \otimes g(B)}VV          @VV{g(A \otimes B)}V   \\
\homa(A) \otimes \homa(B)
@>t_{A,B}>>
\homa(A \otimes B)
\end{CD}
\end{equation*}
}
\end{enumerate}

In \cite{JSGAL} was completely described the Galois group $G^{(0)}_{f}$
of the subcategory $\EE^{(0)}_{f}$ of $\EE^{(0)}$ made up of fuchsian
equations. From \cite{JSFIL}, one could (trivially) deduce the
Galois group $G^{(0)}_{p,1}$ of the category $\EE^{(0)}_{p,1}$ of
pure objects with integral slopes. Here, we will describe the Galois 
group $G^{(0)}$ of $\EE^{(0)}_{1}$. The extension of these results
to the case of non integral slopes should not involve new ideas on
the analytic side, but will have to take in account the results of
van der Put and Reversat in \cite{vdPR}.


\subsection{Galois group and Galois action}
\label{subsection:Galoisaction}

\begin{thm}
\label{theorem:structuregroupe}
The structure of the Galois group $G^{(0)}_{1}$ is as follows:
\begin{eqnarray*}
G^{(0)}_{1} & = & \St \rtimes G_{p_{1}}^{(0)}
\text{~~~(total Galois group with integral slopes)}, \\
G_{p,1}^{(0)} & = & T_{1}^{(0)} \times G_{f}^{(0)}
\text{~~~(pure Galois group with integral slopes)}, \\
T_{1}^{(0)} & = & \C^{*} 
\text{~~~(theta torus with integral slopes)}, \\
G_{f}^{(0)} & = & G_{f,s}^{(0)} \times G_{f,u}^{(0)}
\text{~~~(fuchsian Galois group)}, \\
G_{f,u}^{(0)} & = & \C
\text{~~~(unipotent component of the fuchsian Galois group)}, \\
G_{f,s}^{(0)} & = & Hom_{gr}(\C^{*}/q^{\Z},\C^{*})
\text{~~~(semisimple component of the fuchsian Galois group)}.
\end{eqnarray*}
\end{thm}

The structure and action of the prounipotent Stokes group $\St$ are 
the subject matter of \cite{RS1} and of section \ref{section:wildgroup} 
of the present paper. We shall presently explain the structure and action 
of the pure group $G_{p,1}^{(0)}$. This means that we should associate
to any object $A$ a representation of $G_{p,1}^{(0)}$ in the space
$\homa(A)$; thus, for any $g \in G_{p,1}^{(0)}$ and 
any matrix $A \in GL_{n}(\Ka)$, we should realize $g(A) \in GL_{n}(\C)$. \\

We start from the standard form (\ref{equation:formestandard}). For each 
of the block matrices $A_{i}$, we write:
$$
A_{i} = A_{i,s} A_{i,u}
$$
its multiplicative Dunford decomposition: $A_{i,s}$ is semisimple,
$A_{i,u}$ is unipotent and they commute.

\begin{enumerate}
\item{Let $g = \gamma \in G_{f,s}^{(0)} = Hom_{gr}(\C^{*}/q^{\Z},\C^{*})$.
The latter is here identified with the group of morphisms from the
abstract group $\C^{*}$ to itself that send $q$ to $1$. We let $\gamma$
act on each $A_{i,s}$ through its eigenvalues: if
$A_{i,s} = P \diag(c_{1},\ldots,c_{r}) P^{-1}$, then
$\gamma(A_{i,s}) = P \diag(\gamma(c_{1}),\ldots,\gamma(c_{r})) P^{-1}$
(it does not depend on the choice of a particular diagonalisation).
Then:
$$
g(A) = \begin{pmatrix}
\gamma(A_{1,s})  & \ldots & \ldots & \ldots & \ldots \\
\ldots & \ldots & \ldots  & 0 & \ldots \\
0      & \ldots & \ldots   & \ldots & \ldots \\
\ldots & 0 & \ldots  & \ldots & \ldots \\
0      & \ldots & 0       & \ldots & \gamma(A_{k,s})    
\end{pmatrix},
$$}
\item{Let $g = \lambda \in G_{f,u}^{(0)} = \C$. Since the $A_{i,u}$ 
are unipotent matrices, the $A_{i,u}^{\lambda}$ are well defined and
we put:
$$
g(A) = \begin{pmatrix}
A_{1,u}^{\lambda}  & \ldots & \ldots & \ldots & \ldots \\
\ldots & \ldots & \ldots  & 0 & \ldots \\
0      & \ldots & \ldots   & \ldots & \ldots \\
\ldots & 0 & \ldots  & \ldots & \ldots \\
0      & \ldots & 0       & \ldots & A_{k,u}^{\lambda}    
\end{pmatrix},
$$}
\item{Let $g = t \in T_{1}^{(0)} = \C^{*}$. This \emph{theta torus}
is the analogue here of the \emph{exponential torus} of the classical
differential Galois theory. Then:
$$
g(A) = \begin{pmatrix}
t^{\mu_{1}} I_{r_{1}}  & \ldots & \ldots & \ldots & \ldots \\
\ldots & \ldots & \ldots  & 0 & \ldots \\
0      & \ldots & \ldots   & \ldots & \ldots \\
\ldots & 0 & \ldots  & \ldots & \ldots \\
0      & \ldots & 0       & \ldots & t^{\mu_{k}} I_{r_{k}}    
\end{pmatrix},
$$}
\end{enumerate}

Note that all these depend on $A_{0}$ only. This is because the
category $\EE^{(0)}_{p,1}$ of pure modules with integral slopes
is equivalent to the category of representations of $G^{(0)}_{p,1}$,
so that giving a representation of the latter group is the same as
giving an object in the former category. We leave as an exercise
for the reader the reconstruction of $A_{0}$ from the representation
described above. For further use, we shall now prove two lemmas about 
the action of $G^{(0)}_{p,1}$ on $\homa(A)$. These
lemmas actually express the ``duality'' of $G^{(0)}_{p,1}$ and
$\EE^{(0)}_{p,1}$. 

\begin{lemma}
\label{lemma:vecteurcovariant}
Let $A$ be in normalized standard form (\ref{equation:formestandard}). 
Let $X \in \homa(A)$ be \emph{covariant} under the
action of $G^{(0)}_{p,1}$, that is, for all $g \in G^{(0)}_{p,1}$, 
the vectors $X$ and $g(A) X$ are colinear. Then there exists 
$i \in \{1,\ldots,k\}$ and $\alpha \in \Sp(A_{i})$ such that: 
$A_{0} X = \alpha z^{\mu_{i}} X$.
\end{lemma}
\Pr
First note that the block decomposition of $A_{0}$ (or, equivalently,
the action of the theta torus) entails a splitting of vector spaces: 
$$
\homa(A) = \C^{n} =
\C^{r_{1}} \oplus \cdots \oplus \C^{r_{k}},
$$
each $A_{i}$ acting upon the corresponding $\C^{r_{i}}$. We can
accordingly write $X = (X_{1},\ldots,X_{k})$ (in row form, instead
of column form, for economy of space). Covariance under the action
of $T_{1}^{(0)}$ say that $(t^{\mu_{1}} X_{1},\ldots,t^{\mu_{k}} X_{k})$
and $(X_{1},\ldots,X_{k})$ are colinear for all $t \in \C^{*}$, which
implies that at most one component $X_{i}$ is non trivial. Then,
covariance under the action of $G_{f,u}^{(0)}$ says that $X_{i}$
is fixed by $A_{i,u}$ (since the latter is unipotent). Last, 
covariance under $G_{f,s}^{(0)}$ implies that $X_{i}$ is an 
eigenvector of $A_{i,s}$. Indeed, this comes from the fact that, 
if $\alpha \neq \alpha'$ are eigenvalues of $A_{i}$, then, by the 
normalization condition, $\alpha q^{\Z} \cap \alpha' q^{\Z} = \emptyset$; 
it is then easy to see that there exists 
$\gamma \in Hom_{gr}(\C^{*}/q^{\Z},\C^{*})$ such that
$\gamma(\alpha) \neq \gamma(\alpha')$, so that $X_{i}$ cannot have 
nontrivial components in both eigenspaces of $A_{i}$. The conclusion 
follows.
\hfill $\Box$ 

\begin{lemma}
\label{lemma:vecteurinvariant}
Let $A$ be in normalized standard form (\ref{equation:formestandard}). 
Let $X \in \homa(A)$ be \emph{invariant} under the
action of $G^{(0)}_{p,1}$, that is, for all $g \in G^{(0)}_{p,1}$, 
the vectors $X$ and $g(A) X$ are equal. Then $A_{0} X = X$.
\end{lemma}
\Pr
The proof is similar, with two adaptations. First, equality of
$(t^{\mu_{1}} X_{1},\ldots,t^{\mu_{k}} X_{k})$
and $(X_{1},\ldots,X_{k})$ entails that at most one component $X_{i}$ 
is non trivial and the corresponding slope is $\mu_{i} = 0$; second,
invariance under $G_{f,s}^{(0)}$ implies that at most one component
of $X_{i}$ (in the eigenspace decomposition) is non trivial, that the
corresponding $\alpha \in \Sp(A_{i})$ is in the kernel of all elements
of $Hom_{gr}(\C^{*}/q^{\Z},\C^{*})$, so it is in $q^{\Z}$, so equal
to $1$ by the normalisation condition.
\hfill $\Box$ \\

Again because of the duality of $G^{(0)}_{p,1}$ and $\EE^{(0)}_{p,1}$,
the conclusions of these two lemmas have useful interpretations. 
The conclusion of lemma \ref{lemma:vecteurcovariant} says that 
the column matrix $X \in M_{n,1}(\C)$ is a morphism from the rank 
one object $(\alpha z^{\mu_{i}}) \in GL_{1}(\Ka)$ into $A_{0}$. 
The conclusion of lemma \ref{lemma:vecteurinvariant} says that 
the column matrix $X \in M_{n,1}(\C)$ is a morphism from the unit 
object $\underline{1} = (1) \in GL_{1}(\Ka)$ into $A_{0}$, \emph{i.e.}
a section $X \in \Gamma(A_{0})$.



\section{The wild fundamental group}
\label{section:wildgroup}


\subsection{The action of the Stokes group}
\label{subsection:Stokesaction}

An element $s \in \St$ is characterized by the following properties:
\begin{enumerate}
\item{To each $A$ in standard form (\ref{equation:formestandard}),
it associates a matrix $s(A) \in \G(\C)$; recall that $\G$ was 
described as the algebraic group of matrices of shape as in equation 
(\ref{equation:automorphisme}).}
\item{If $A = A_{0}$, that is, if $A$ is pure, then $s(A) = I_{n}$.}
\item{Functoriality and tensor compatibility are defined as in section 
\ref{subsection:fiberfunctor}.}
\end{enumerate}

Since $\St$ is a prounipotent proalgebraic group, it is convenient
to study it through its Lie algebra $\st$. (The underlying formalism
is expounded in the appendix of \cite{DG}.) An element $D \in \st$
is also a natural transformation of $\homa$. 
It associates to each object $A$ an endomorphism
$D(A) \in \mathcal{L}(\homa(A)) = M_{n}(\C)$,
subject to the following conditions:
\begin{enumerate}
\item{For each $A$ in standard form (\ref{equation:formestandard}),
the matrix $D(A) \in M_{n}(\C)$ is in $\g(\C)$; recall that $\g$ was
described as the Lie algebra of matrices of shape as in equation 
(\ref{equation:endomorphisme}).}
\item{If $A = A_{0}$, that is, if $A$ is pure, then $D(A) = 0_{n}$.}
\item{Functoriality is defined as in section \ref{subsection:fiberfunctor}.}
\item{Tensor compatibility is that of ``Lie-like elements'' (as in
\cite{Simpson}, \S 6): for any $A \in GL_{n}(\Ka)$ and $B \in GL_{p}(\Ka)$,
one should have, up to natural identifications:
$D(A \otimes B) = D(A) \otimes I_{p} + I_{n} \otimes D(B)$.
Thus, $D$ behaves like a derivation.}
\end{enumerate}

In \cite{RS1}, we have produced many elements of $\St$ and of $\st$.
However, for a given basepoint $a \in \C^{*}$, these do not operate
on the whole of $\EE^{(0)}_{1}$ but on a tannakian subcategory of it.
Therefore, the way of using them is the following: given an equation
$A$ of interest, proposition 4.2 of \emph{loc. cit.} yields an explicit 
criterion to select adequate basepoints (these are generically adequate).
Then all the constructions that follow make sense in the tannakian
subcategory of $\EE^{(0)}_{1}$ generated by $A$. This means that
each time we shall evaluate a meromorphic function at $a$, this
will be possible. Henceforth, we shall not anymore discuss this
matter. \emph{We assume that the basepoint has been chosen so
that all the objects we deal with are compatible with it.} \\

In \cite{JSSTO} and \cite{RS1}, we defined an explicit finite subset 
$\Sigma_{A_{0}}$ of $\Eq$ and proved:

\begin{thm}
\label{thm:ScF}
Let $\overline{c} \in \Eq \setminus \Sigma_{A_{0}}$. Then, there is 
a unique $F: A_{0} \rightarrow A$ such that $F \in \G(\Kw)$, with
poles only on $[-c;q] = -c q^{\Z}$ and such that, for $1 \leq i < j \leq k$, 
the poles of $F_{i,j}$ have multiplicity $\leq \mu_{j} - \mu_{i}$.
\end{thm}

We write this meromorphic isomorphism $S_{\overline{c}} \hat{F}_{A}$
and see it as some kind of \emph{summation of $\hat{F}_{A}$ in the 
direction $\overline{c} \in \Eq$}. Therefore, changing direction
of summation, we may define, for every
$\overline{c},\overline{d} \in \Eq \setminus \Sigma_{A_{0}}$:
$$
S_{\overline{c},\overline{d}} \hat{F}_{A} =
(S_{\overline{c}} \hat{F}_{A})^{-1} S_{\overline{d}} \hat{F}_{A},
$$
some kind of ``ambiguity of summation'', that is, a Stokes operator.
It is plainly a meromorphic automorphism of $A_{0}$. We also proved
in \emph{loc. cit.}:

\begin{prop}
If moreover $\overline{a} \neq \overline{c},\overline{d}$, then
$A \leadsto S_{\overline{c},\overline{d}} \hat{F}_{A}(a)$ is
an element of $\St$. In particular,
$S_{\overline{c},\overline{d}} \hat{F}_{A}(a) \in \St(A)$. 
(Recall that we implicitly restrict ourselves to a subcategory 
of $\EE^{(0)}_{1}$ where everything is defined.)
\end{prop}

For the following corollary, we fix an arbitrary direction of summation
$\overline{c_{0}} \in \Eq$, again to be considered as a choice of
basepoint (and inessential).

\begin{cor}
Putting
$\Lsca(A) =
\log (S_{\overline{c_{0}},\overline{c}} \hat{F}_{A}(a)) \in \st(A)$
yields a family of elements of elements of $\st(A)$. Moreover,
$A \leadsto \Lsca(A)$ is an element of $\st$. 
(We omit $\overline{c_{0}}$ in the notation.)
\end{cor}

The above family is a meromorphic map from $\Eq$ to a vector space,
hence one can take residues. Define the $q$-alien derivations by the 
formula:
$$
\Derc(A) = Res_{\overline{d} = \overline{c}} LS_{\overline{d},a}(A).
$$
(We do not mention the arbitrary basepoints $\overline{c_{0}}, a$ in 
the notation.)
Of course, for $\overline{c} \not\in \Sigma_{A_{0}}$, we have 
$\Derc(A) = 0$. Another result we need from \cite{RS1} is:

\begin{thm}
One has $\Derc(A) \in \st(A)$. More precisely, $A \leadsto \Derc(A)$ 
is an element of $\st$.
\end{thm}

Since $\St$ is a normal subgroup of $G^{(0)}_{1}$, it admits a
conjugation action by $G^{(0)}_{p,1}$. This can be transferred
to the Lie algebra $\st$. Because of the action by the theta torus
$T^{(0)}_{1} = \C^{*}$, we thus have a spectral decomposition:
$$
\st = \bigoplus_{\delta \geq 1} \st^{\delta},
$$
and each alien derivation admits a canonical decomposition:
$$
\Derc = \bigoplus_{\delta \geq 1} \Derdc,
$$
where $\Derdc(A) \in \st^{\delta}(A)$ has only non null blocks for
$\mu_{j} - \mu_{i} = \delta$. Each $t \in T^{(0)}_{1}$ acts on
$\st^{\delta}$ by multiplication by $t^{\delta}$, and carries 
$\Derdc(A)$ to $t^{\delta} \Derdc(A)$. \\

\Rems
\begin{enumerate}
\item{The theta torus actually operates on each $\homa(A) = \homa(A_{0})$ 
and, being semi-simple, splits it into the direct sum of its eigenspaces:
one for each slope $\mu$, with rank $r(\mu)$. The corresponding increasing
filtration comes from the filtration by the slopes:
$$
\homa(A)_{\geq \mu} = \homa(A_{\geq \mu}).
$$
The elements of the group $\G(\C)$ are the automorphisms of $\homa(A)$
which respect that filtration and are trivial (\emph{i.e.} the identity)
on the associated graded space. The elements of the algebra $\G(\C)$ are 
the endomorphisms of $\homa(A)$ which respect that filtration and are 
trivial (\emph{i.e.} null) on the associated graded space.}
\item{From this, we deduce a spectral decomposition:
$$
\g = \bigoplus_{\delta \geq 1} \g^{\delta},
$$
from which the decomposition of $\st$ follows.}
\item{Putting 
$\g^{\geq \delta} = \sum\limits_{\delta' \geq \delta} \g^{\delta'}$
defines a filtration of the Lie algebra $\g$ by ideals. 
Putting $\G^{\geq \delta} = I_{n} + \g^{\geq \delta} = \exp \g^{\geq \delta}$ 
then defines a filtration of $\G$ by normal subgroups.}
\item{Similarly, we can decompose each $(i,j)$ block of $\Derdc(A)$
(where $\mu_{j} - \mu_{i} = \delta$) into subblocks indexed
by pairs $(\alpha,\beta) \in \Sp(A_{i}) \times \Sp(A_{j})$.
Each $\gamma \in G^{(0)}_{f,s}$ then multiplies the corresponding
subblock of $\Derdc(A)$ by $\dfrac{\gamma(\alpha)}{\gamma(\beta)} \cdot$
This may be deduced as above from the action of $G^{(0)}_{f,s}$ on
$\homa(A_{0})$ and a corresponding splitting of each of the 
$\g^{\delta}$.}
\end{enumerate}


\subsection{The density theorem}
\label{subsection:density}


\subsubsection{Plain density theorem}

The wild monodromy group actually is the Lie subalgebra of $\st$
generated by the $q$-alien derivations $\Derc$. The justification
of the name is that its definition has a transcendental character 
and the following result.

\begin{thm}[Density theorem]
\label{thm:density}
(i) The subgroup of $\St$ associated with the wild monodromy group 
(as defined above), together with the pure group $G^{(0)}_{p,1}$, 
generate a Zariski-dense subgroup of the whole Galois group $G^{(0)}_{1}$. \\
(ii) The $\Derc$ together with all their conjugates under the action 
of $G^{(0)}_{p,1}$ generate a Zariski-dense Lie subalgebra of $\st$. 
\end{thm}
\Pr
Actually, (i) is but a rephrasing of (ii) and we shall prove the latter.
We shall use Chevalley's criterion in the following form: \\
{\sl For a subset $H \subset G^{(0)}_{1}$ to generate a Zariski-dense 
subgroup of $G^{(0)}_{1}$, it is sufficient that, for each object
$A$ and each line $D \subset \homa(A)$ which is
invariant under the action of all elements of $H$, then $D$ is
actually invariant under the action of $G^{(0)}_{1}$.} Our way of 
using it is similar to that in \cite{JSGAL} (2.2.3.3 and 3.1.2.3). \\

We take
$H = G^{(0)}_{p,1} \; \cup \; \exp\left(\{\Derc \mid c \in \C^{*}\}\right)$.
If we choose a generator $X$ of the line $D$, the assumption is that
$X$ is covariant under $G^{(0)}_{p,1}(A) = G^{(0)}_{p,1}(A_{0})$ on 
the one hand, under all the $\Derc(A)$ on the other hand. Since the
latter are nilpotent, this means that all $\Derc(A) X = 0$. Then,
must prove that for all $D \in \st$, one has $D(A) X = 0$. \\

Using lemma \ref{lemma:vecteurcovariant}, along with its proof and 
notations, we may write (in row form) $X = (0,\ldots,X_{i},\ldots,0)$, 
where the components have sizes $r_{1},\ldots,r_{k}$ and where 
$A_{i} X_{i} = c X_{i}$ for some $c \in \Sp(A_{i})$, so that
$A_{0} X = c z^{\mu_{i}} X$.  \\

Now, we note that components of slopes $> \mu_{i}$ are neither involved
in the assumptions nor in the conclusion, so that one may as well assume 
from start that $i = k$. Indeed, write $n' = r_{1} + \cdots + r_{i}$ the 
size of the components corresponding to slopes $\leq \mu_{i}$ (equivalently, 
the rank of the submodule $M'_{\leq \mu_{i}} \subset M = M_{A}$ of slopes 
$\leq \mu_{i}$ in the slope filtration), $A'$ the corresponding submatrix 
of $A$ (so that $M' = M_{A'}$) and $X' = (0,\ldots,X_{i})$ the corresponding
subvector of $X$. The matrix 
$\Phi = \begin{pmatrix} I_{n'} \\ 0 \end{pmatrix} \in M_{n,n'}(\C)$
is a morphism $\Phi: A' \rightarrow A$ (it is the inclusion
$M' \subset M$) and $\Phi X' = X$. For all $g \in G^{(0)}_{1}$,
one has (by functoriality) $g(A) \Phi = \Phi g(A')$ (here, one
has $\Phi(a) = \Phi$). The reader
will check that $A'$ and $X'$ satisfy the same assumption as $A$
and $X$, and that it is enough to prove the conclusions for them. \\

So we assume from now on that $X = (0,\ldots,X_{k})$, that 
$A_{k} X_{k} = c X_{k}$ for some $c \in \Sp(A_{k})$, so that
$A_{0} X = c z^{\mu_{k}} X$. Then $X: (c z^{\mu_{k}}) \rightarrow A_{0}$
is an analytic morphism, and therefore 
$G = \hat{F}_{A} X: (c z^{\mu_{k}}) \rightarrow A$ is a formal morphism.
We shall prove below (lemma \ref{lemma:hatFXanalytique}) that it is 
actually an analytic morphism. Therefore, taking $D \in \st$ and
using functoriality, we get the commutative diagram:
\begin{equation*}
\begin{CD}
\C                      @>G_{0}(a)>>            \homa(A) \\
@V{D(c z^{\mu_{k}})}VV                          @VV{D(A)}V      \\
\C                      @>G_{0}(a)>>            \homa(A)
\end{CD}
\end{equation*}
The matrix $G_{0}$ is the graded part of the column matrix $\hat{F}_{A} X$,
that is, in row notation, $X$ itself. Since $D \in \st$ and since the source 
object $(c z^{\mu_{k}})$ is pure, one has $D(c z^{\mu_{k}}) = 0$. Hence
we get $D(A) X = 0$ as wanted.
\hfill $\Box$

\begin{lemma}
\label{lemma:hatFXanalytique}
The matrix $\hat{F}_{A} X$ is analytic and the summations
$S_{\overline{c}} \hat{F}_{A}$ do not depend on the direction 
$\overline{c} \in \Eq$, and they are all equal to $\hat{F}_{A} X$
(that is, its ``classical summation'', as a \emph{convergent}
power series).
\end{lemma}
\Pr
First fix a direction
$\overline{c} \in \Eq$ and write $F = S_{\overline{c}} \hat{F}_{A}$.
Likewise, write $G = \hat{F}_{A}$ for short. The components $F_{i,j}$
and $G_{i,j}$ satisfy equations (\ref{equation:calculrecursifhat(F)}). 
We are interested in the $F_{i,k} X_{k}$ and the $G_{i,k} X_{k}$ for 
$1 \leq i < k$ (for $i = k$, both equal $X_{k}$). \\

We shall put: $Y_{i} = F_{i,k} X_{k}$, $Z_{i} = G_{i,k} X_{k}$ 
and $V_{i} = U_{i,k} X_{k}$. Then, for $1 \leq i < k$, multiplying 
(\ref{equation:calculrecursifhat(F)}) by $X_{k}$ on the right and 
taking in account the equalities $\sq X_{k} = X_{k}$ and 
$A_{k} X_{k} = c X_{k}$, one gets:
$$
c z^{\mu_{k}} (\sq Y_{i}) - z^{\mu_{i}} A_{i} Y_{i} =
\sum_{i < j < k} U_{i,j} Y_{j} + V_{i} \quad \text{and} \quad
c z^{\mu_{k}} (\sq Z_{i}) - z^{\mu_{i}} A_{i} Z_{i} =
\sum_{i < j < k} U_{i,j} Z_{j} + V_{i}.
$$

On the other hand, we shall have to use the assumptions: $\Derc(A) X = 0$. 
Since $X = (0,\ldots,X_{k})$, this means that, for each $i < k$, one has
$(\Derc(A))_{i,k} X_{k} = 0$. Writing for short
$F_{d} = S_{\overline{d}} \hat{F}_{A}(a)$,
$F_{0} = S_{\overline{c_{0}}} \hat{F}_{A}(a)$ and 
$L = \Lsda(A) =
\log (S_{\overline{c_{0}},\overline{d}} \hat{F}_{A}(a))$
(see the notations of section \ref{subsection:Stokesaction}),
we shall see in lemma \ref{lemma:CHfaible} below that, for $i < k$:
$$
L_{i,k} = (F_{d})_{i,k} - (F_{0})_{i,k} +
\sum_{i < j < k} M_{i,j,k} \bigl((F_{d})_{j,k} - (F_{0})_{j,k}\bigr),
$$
where the $M_{i,j,k}$ are some arbitrary matrices (their values are
inessential here). \\

We use a downward induction on $i$. For $i = k-1$:
$$
c z^{\mu_{k}} (\sq Y_{k-1}) - z^{\mu_{k-1}} A_{k-1} Y_{k-1} = V_{k-1}
\quad \text{and} \quad
c z^{\mu_{k}} (\sq Z_{k-1}) - z^{\mu_{k-1}} A_{k-1} Z_{k-1} = V_{k-1}.
$$
Thus, $Z_{k-1}$ is the formal solution and $Y_{k-1}$ the solution summed 
in direction $\overline{c}$ of the equation 
$c z^{\mu_{k}} (\sq Y) - z^{\mu_{k-1}} A_{k-1} Y = V_{k-1}$. 
On the other hand, from the formula above, one has:
$L_{k-1,k} = (F_{d})_{k-1,k} - (F_{0})_{k-1,k}$.
Taking residue at $\overline{d} = \overline{c}$ and multiplying
at right by the constant vector $X_{k}$, one gets:
$$
(\Derc(A))_{k-1,k} X_{k} = 
Res_{\overline{d} = \overline{c}} LS_{k-1,k} X_{k} =
Res_{\overline{d} = \overline{c}} (F_{d})_{k-1,k} X_{k} =
Res_{\overline{d} = \overline{c}} Y_{k-1}.
$$
This means that the residues of resummed solutions of the equation
just mentioned are all $0$. According to the results of section 4.2
of \cite{RS1}, this implies that $Y_{k-1}$ is analytic near $0$ (it
has no poles other than $0$), that it does not depend on the direction
of summation $\overline{d}$, and that it is equal to $Z_{k-1}$. This
completes the first step of the induction. \\

Now take $i < k-1$ and assume that the property has been proved
for all $j > i$. Consider the equation of which $Y_{i}$ is solution.
Its right hand member $\sum\limits_{i < j < k} U_{i,j} Y_{j} + V_{i}$
is analytic, by the induction hypothesis (analyticity part). The
residue at $\overline{d} = \overline{c}$ of $Y_{i}$ is equal to,
$Res_{\overline{d} = \overline{c}} (F_{d})_{i,k} X_{k}$, thus to:
$$
Res_{\overline{d} = \overline{c}} LS_{i,k} X_{k} = 
(\Derc(A))_{i,k} X_{k} = 0.
$$ 
This is because all other terms in the formula taken from lemma 
\ref{lemma:CHfaible} have at right a factor 
$\bigl((F_{d})_{j,k} - (F_{0})_{j,k}\bigr) X_{k} = Y_{j} - Y_{j}$,
since $Y_{j}$ does not depend on the direction of summation. Thus
we have again a solution $Y_{i}$ with all residues null, so that it 
is analytic and independent of the direction of summation by
\emph{loc. cit.}.
\hfill $\Box$

\begin{lemma}
\label{lemma:CHfaible}
With the notations of section \ref{subsection:Stokesaction},
one has, for $i < k$:
$$
\log (S_{\overline{c_{0}},\overline{d}} \hat{F}_{A}(a))_{i,k} = 
(S_{\overline{d}} \hat{F}_{A}(a))_{i,k} - 
(S_{\overline{c_{0}}} \hat{F}_{A}(a))_{i,k} +
\sum_{i < j < k} M_{i,j,k} \bigl(
(S_{\overline{d}} \hat{F}_{A}(a))_{j,k} - 
(S_{\overline{c_{0}}} \hat{F}_{A}(a))_{j,k}
\bigr),
$$
where the $M_{i,j,k}$ are some arbitrary matrices.
\end{lemma}
\Pr
We write $A = S_{\overline{d}} \hat{F}_{A}(a)$,
$B = S_{\overline{c_{0}}} \hat{F}_{A}(a)$ and $C = A - B$, 
which is strictly upper triangular by blocks. Then:
$$
\log(B^{-1} A) = \log(I_{n} + B^{-1} C) = 
\sum_{p \geq 1} \dfrac{(-1)^{p-1}}{p} (B^{-1} C)^{p},
$$
from which the equality of blocks:
$$
\bigl(\log(B^{-1} A)\bigr)_{i,k} = C_{i,k} +
\sum_{i < j < k} M_{i,j,k} C_{j,k}
$$
follows easily.
\hfill $\Box$


\subsubsection{Functorial density theorem}

In section \ref{subsection:freeness}, we shall describe how the Zariski 
generators $\Derc$ of $\st$ (theorem \ref{thm:density}) are related. 
For that, we shall first give a more functorial version of the density 
theorem. \\

Since $\EE^{(0)}_{p,1}$ and $\EE^{(0)}_{1}$ are respectively isomorphic 
to the category of (finite dimensional complex) representations of 
$G^{(0)}_{p,1}$ and $G^{(0)}_{1} = \St \rtimes G^{(0)}_{p,1}$, and 
since finite dimensional representations of the prounipotent
proalgebraic group $\St$ are equivalent to finite dimensional 
representations of the pronilpotent proalgebraic Lie algebra $\st$,
the tannakian category $\EE^{(0)}_{1}$ admits an alternative ``mixed''
description, which runs as follows:
\begin{enumerate}
\item{Objects are pairs $A = \left(A_{0},(D(A))_{D \in \st}\right)$,
where $A_{0}$ is some object of $\EE^{(0)}_{p,1}$, \emph{e.g.} a
matrix in pure standard form (\ref{equation:formestandardpure}),
and where each $D(A) \in \g(\C)$.}
\item{Morphisms from $A = \left(A_{0},(D(A))_{D \in \st}\right)$
to $B = \left(B_{0},(D(B))_{D \in \st}\right)$ are morphisms
$F_{0}: A_{0} \rightarrow B_{0}$ in $\EE^{(0)}_{p,1}$ such that,
for each $D \in \st$, one has $D(B) F_{0}(a) = F_{0}(a) D(A)$.
(Recall that an arbitrary basepoint $a \in \C^{*}$ has been chosen
once for all.)}
\item{The tensor product of $A = \left(A_{0},(D(A))_{D \in \st}\right)$
and $B = \left(B_{0},(D(B))_{D \in \st}\right)$ is the object
$C = \left(C_{0},(D(C))_{D \in \st}\right)$, with the previous rule
$C_{0} = A_{0} \otimes B_{0}$ from $\EE^{(0)}_{p,1}$, and with the
``Lie-like element'' rule 
$D(A \otimes B) = D(A) \otimes I_{p} + I_{n} \otimes D(B)$.
The unit is $\underline{1} = \left((1),(0)_{D \in \st}\right)$.
The dual of $A$ is
$A^{\vee} = \left(A_{0}^{\vee},(- {}^{t}D(A))_{D \in \st}\right)$.
The space of sections of $A$ is
$\Gamma(A) = Hom(\underline{1},A) = 
\{X_{0} \in \Gamma(A_{0}) \mid \forall D \in \st \;,\; D(A) X_{0} = 0\}$.
(Recall that 
$\Gamma(A_{0}) = \{X_{0} \in \Ka^{n} \mid \sq X_{0} = A_{0} X_{0}\}$.)}
\item{There is a a fiber functor $\homa(A) \underset{def}{=} \homa(A_{0})$.}
\end{enumerate}
To be complete, such a description should take in account the adjoint 
action of $G^{(0)}_{p,1}$ on $\st$, which is, for all $A$, the restriction
of the action of $G^{(0)}_{p,1}$ on $\g(\C)$. For instance, from the
action of the theta torus, one draws the graduation 
$\st = \bigoplus\limits_{\delta \geq 1} \st^{\delta}$, whence 
decompositions $D(A) = \sum\limits_{\delta \geq 1} D^{\delta}(A)$,
where each $D^{\delta}(A) \in \g^{\delta}(\C)$. We shall take
in account the adjoint action of $G_{f}^{(0)}$ later. \\

We would like to consider the $\Derdc(A)$ as encoding a Lie algebra 
representation from the free Lie algebra $L$ generated by the family 
of symbols $(\Derdc)_{\delta \geq 1,\overline{c} \in \Eq}$, and so
describe $\EE^{(0)}_{1}$ as the category of representations of
$L \rtimes G^{(0)}_{p,1}$ in a way similar to that above. This would
require some other tools (see the conclusion of the paper). As a
substitute, we define a new tannakian category $\EE'$ as follows:
\begin{enumerate}
\item{Objects are pairs:
$$
A = \left(A_{0},(\Derdc(A))_{\delta \geq 1,\overline{c} \in \Eq}\right),
$$
where $A_{0}$ is in pure standard form (\ref{equation:formestandardpure}),
and where each $\Derdc(A) \in \g^{\delta}(\C)$.}
\item{Morphisms from $A$
to $B = \left(B_{0},(\Derdc(B))_{\delta \geq 1,\overline{c} \in \Eq}\right)$ 
are morphisms
$F_{0}: A_{0} \rightarrow B_{0}$ in $\EE^{(0)}_{p,1}$ such that,
for each $\delta \geq 1,\overline{c} \in \Eq$, one has 
$\Derdc(B) F_{0}(a) = F_{0}(a)\Derdc(A)$.}
\item{The tensor product of $A$ and $B$ is the object
$C = \left(C_{0},(\Derdc(B))_{\delta \geq 1,\overline{c} \in \Eq}\right)$, 
with $C_{0} = A_{0} \otimes B_{0}$ and
$\Derdc(A \otimes B) = \Derdc(A) \otimes I_{p} + I_{n} \otimes \Derdc(B)$.
The unit and dual are described as before. The space of sections of $A$ is
$\Gamma(A) = Hom(\underline{1},A) = \{X_{0} \in \Gamma(A_{0}) \mid 
\forall \delta \geq 1,\overline{c} \in \Eq \;,\; \Derdc(A) X_{0} = 0\}$.}
\item{There is a a fiber functor $\homa(A) \underset{def}{=} \homa(A_{0})$.}
\end{enumerate}
For the time being, we do not take in account the action of $G_{f}^{(0)}$. \\

We now consider the functor 
$A \leadsto \F(A) = 
\left(A_{0},(\Derdc(A))_{\delta \geq 1,\overline{c} \in \Eq}\right)$
from $\EE^{(0)}_{1}$ to $\EE'$. It is plainly an exact faithful
$\otimes$-functor.

\begin{thm}[Functorial density theorem]
\label{thm:functorialdensity}
The functor $\F$ is fully faithful.
\end{thm}
\Pr
To prove that $Hom(A,B) \rightarrow Hom(\F(A),\F(B))$ is onto,
we draw on the identifications $Hom(A,B) = \Gamma(A^{\vee} \otimes B)$
and $Hom(\F(A),\F(B)) = \Gamma(\F(A)^{\vee} \otimes \F(B))$. Since
$\F$ is a $\otimes$-functor, the latter is identified with
$\Gamma(\F(A^{\vee} \otimes B))$, so that we are left check that, 
for any $A$, the map $\Gamma(A) \rightarrow \Gamma(\F(A))$ is onto. \\

That map sends a vector
$X \in \Ka^{n}$ such that $\sq X = A X$ to its graded part
$X_{0} \in \Ka^{n}$. The vector $X_{0}$ has the same null slope component
as $X$ and is zero elsewhere. It satisfies $\sq X_{0} = A_{0} X_{0}$
and $\forall \delta \geq 1,\overline{c} \in \Eq \;,\; \Derdc(A) X_{0} = 0$.
If we start from such a vector $X_{0}$, lemma \ref{lemma:vecteurinvariant}
tells us that it comes indeed from some $X \in \Gamma(A)$.
\hfill $\Box$


\subsection{A freeness theorem}
\label{subsection:freeness}

We now shall describe the (essential) image of the functor $\F$,
or, what amounts to the same, which families 
$(\Derdc(A))_{\delta \geq 1,\overline{c} \in \Eq}$
can be realized for a given $A_{0}$ in $\EE^{(0)}_{p,1}$.
To understand what is going on, we start with the first level,
which is easier.


\subsubsection{The first level}

In theorem \ref{thm:ScF}, we obtained the $F_{i,j}$ blocks of
$S_{\overline{c}} \hat{F}_{A}$ as solutions of the following
equations:
$$
(\sq F_{i,j}) z^{\mu_{j}} A_{j} - z^{\mu_{i}} A_{i} F_{i,j} =
\sum_{i < l < j} U_{i,l} F_{l,j} + U_{i,j}.
$$
We consider the first non trivial level in the computation
of $S_{\overline{c}} \hat{F}_{A}$, that is:
$\delta_{0} = \min\limits_{i < j}(\mu_{j} - \mu_{i})$. 
For a block $F_{i,j}$ of level $\mu_{j} - \mu_{i} = \delta_{0}$,
there is no $F_{l,j} \neq 0$ for $i < l < j$, so that the second hand 
member in the equation above is $U_{i,j}$, which is analytic near $0$. 
In that case, there is a solution $F_{i,j}$ with poles on $[-c;q]$ and 
multiplicity $\leq \delta_{0}$ for any $\overline{c} \in \Eq$ which satisfies 
the non-resonancy condition:
$$
\forall \alpha \in \Sp(A_{i}) \;,\; \forall \beta \in \Sp(A_{j}) \;,\;
\alpha c^{\mu_{i}} \not\equiv \beta c^{\mu_{j}} \pmod{q^{\Z}}. 
$$
We recall briefly, from \cite{RS1}, how this was computed. One puts 
$F_{i,j} = \dfrac{G_{i,j}}{\theta_{c}^{\delta_{0}}}$, (the function
$\theta_{c}$ has been defined in section \ref{subsection:generalnotations}).
We thus look for $G_{i,j}$ holomorphic on $\C^{*}$ and satisfying:
$$
c^{\delta_{0}} (\sq G_{i,j}) A_{j} - A_{i} G_{i,j} =
z^{-\mu_{i}} U_{i,j} \theta_{c}^{\delta_{0}} = \sum_{n \in \Z} v_{n} z^{n}.
$$
Writing the Laurent series $G_{i,j} = \sum g_{n} z^{n}$, we are
left to solve, for each $n \in \Z$:
$$
c^{\delta_{0}} q^{n} g_{n} A_{j} - A_{i} g_{n} = 
v_{n} \in \Ma_{r_{i},r_{j}}(\C).
$$
If $\Sp(c^{\delta_{0}} q^{n} A_{j}) \cap \Sp(A_{i}) = \emptyset$,
which is is just the non-resonancy condition above, then, for each $n$, 
this admits a unique solution . \\

Using the notations given at the end of \ref{subsection:Stokesaction},
we see that, for any $\overline{d} \not\in \Sigma_{A_{0}}$, one has
$S_{\overline{d}} \hat{F}_{A}(a) \in \G^{\geq \delta_{0}}(\C)$.
Provisionally call $f^{\delta}_{\overline{d}}$ its component at
level $\delta$. A small computation shows that $LS_{\overline{d},a}(A)$ 
is in $\g^{\geq \delta_{0}}(\C)$ and that its component at level $\delta_{0}$
is $f^{\delta_{0}}_{\overline{d}} - f^{\delta}_{\overline{c_{0}}}$. Thus,
$\Der^{(\delta_{0})}_{\overline{c}}(A) = 
Res_{\overline{d} = \overline{c}} f^{\delta_{0}}_{\overline{d}}$. \\

>From the previous computation, we now conclude that, for 
$\mu_{j} - \mu_{i} = \delta_{0}$, the $(i,j)$ block of 
$\Der^{(\delta_{0})}_{\overline{c}}(A)$ is trivial for
non-resonant directions, \emph{i.e.} if
$\Sp(c^{\delta_{0}} q^{n} A_{j}) \cap \Sp(A_{i}) = \emptyset$.
This is the necessary conditions we were looking for. It is not hard 
to see (and it will come as a particular case of the following sections)
that these are indeed the only conditions on the first level.


\subsubsection{Structure of the $q$-alien derivations at an arbitrary level}
\label{subsubsection:arbitrarylevel}

We are led to introduce some more notations. We first refine 
the spectral decomposition of $\homa(A)$ under the action of 
the theta torus by taking in account the action of $G^{(0)}_{f,s}$,
the semi-simple component of the fuchsian group. From the equivalence: 
$$
\alpha \equiv \beta \pmod{q^{\Z}} \Longleftrightarrow
\forall \gamma \in G^{(0)}_{f,s} \;,\; \gamma(\alpha) = \gamma(\beta),
$$
we see that the action of $G^{(0)}_{f,s}$ splits each eigenspace 
under $\C^{*}$ corresponding to the slope $\mu_{i}$ into a sum 
indexed by the $\overline{\alpha} \in \overline{\Sp(A_{i})}$. 
Precisely, if $V = \homa(A)$, then one may write:
$$
V = \bigoplus V^{(\mu)},
$$
where $\mu$ runs through the set of slopes of $A$, and, for each $\mu$:
$$
V^{(\mu)} = \bigoplus V^{(\mu,\overline{\alpha})},
$$
where $\alpha$ runs through $\Sp(A_{i})$ if $\mu = \mu_{i}$
in our usual notations. \\

To be able to carry this splitting to matrices, we fix an arbitrary
linear order on $\Eq$ and assume the order on indices is compatible
with that arbitrary order. The corresponding adjoint action of 
$G^{(0)}_{f,s}$ on $\g(A)$ then allows one to define the eigenspaces: 
$$
\g^{(\delta,\overline{c})}(\C) =
\{M \in \g(\C) \mid M 
\text{~is trivial out of the~}
(\mu_{i},\overline{\alpha},\mu_{j},\overline{\beta}) 
\text{~components such that~}
\alpha c^{\mu_{i}} \equiv \beta c^{\mu_{j}} \pmod{q^{\Z}}
\}.
$$
This can be non-trivial only if $\overline{c} \in \Sigma_{A_{0}}^{\delta}$,
where:
$$
\Sigma_{A_{0}}^{\delta} = \{\overline{c} \in \Eq \mid
\exists i < j \text{~such that~} \mu_{j} - \mu_{i} = \delta \text{~and~} 
\dfrac{\overline{\alpha}}{\overline{\beta}} = \overline{c^{\delta}}\}.
$$
By definition, 
$\Sigma_{A_{0}} = \bigcup\limits_{\delta \geq 1} \Sigma_{A_{0}}^{\delta}$.
then:
$$
\g^{\delta}(\C) = 
\bigoplus_{\overline{c} \in \Sigma_{A_{0}}^{\delta}}
\g^{(\delta,\overline{c})}(\C).
$$

\Ex
>From the previous paragraph, it follows that on the first non trivial
level, $\Derdc(A) \in \g^{\delta,\overline{c}}(\C)$. The difficulty
is to properly generalize this fact to upper levels. \\

\Rem
\label{remark:espacefixe}
The equality 
$\overline{\alpha} \overline{c}^{\mu_{i}} = 
\overline{\beta} \overline{c}^{\mu_{j}}$
is equivalent to: 
$\forall \gamma \in G^{(0)}_{f,s} \;,\; 
\gamma(\alpha c^{\mu_{i}}) = \gamma(\beta c^{\mu_{j}})$.
Thus, $\g^{\delta,\overline{c}}(\C)$ can be characterized as
the common fixed space of all the
$\left(\gamma(\overline{c})^{-1},\gamma,0\right) \in G^{(0)}_{p,1}$,
where $\gamma$ runs through $G^{(0)}_{f,s}$. \\

Now let $A, A'$ be matrices in standard form with the same graded part
$A_{0}$. From \cite{JSSTO} and \cite{RS1}, we have the following 
generalisation of theorem \ref{thm:ScF}: for each 
$\overline{c} \in \Eq \setminus \Sigma_{A_{0}}$, there exists a unique
meromorphic morphism $F:A \rightarrow A'$ in $\G(\Kw)$, with poles
on $[-c;q]$ and with multiplicities prescribed as in the theorem.
We write it $S_{\overline{d}} \hat{F}_{A,A'}$. One then has:
$$
S_{\overline{d}} \hat{F}_{A,A'} = 
S_{\overline{d}} \hat{F}_{A'}
\left(S_{\overline{d}} \hat{F}_{A}\right)^{-1}.
$$
Assume now that $A \equiv A' \pmod{\g^{\geq \delta}(\Ka)}$, that is, 
$A$ and $A'$ have the same over-diagonals at levels $< \delta$.
The components $F_{i,j}$ of $S_{\overline{d}} \hat{F}_{A,A'}$ 
for $0 < \mu_{j} - \mu_{i} < \delta$ are solutions of the equations:
$$
(\sq F_{i,j}) z^{\mu_{j}} A_{j} - z^{\mu_{i}} A_{i} F_{i,j} = 0.
$$
Therefore, they are null (\emph{cf. loc. cit.}). This implies:
$$
S_{\overline{d}} \hat{F}_{A,A'} \in \g^{\geq \delta}(\Kw).
$$
>From the equality: $S_{\overline{d}} \hat{F}_{A'} = 
S_{\overline{d}} \hat{F}_{A,A'} S_{\overline{d}} \hat{F}_{A}$,
we deduce:
$$
S_{\overline{d}} \hat{F}_{A} \equiv S_{\overline{d}} \hat{F}_{A'} 
\pmod{\g^{\geq \delta}(\Kw)}.
$$ 

\begin{prop}
Let $f_{A,A',\overline{d}}$ be the component at level $\delta$
of $S_{\overline{d}} \hat{F}_{A,A'}(a)$. Then:
$$
\Derdc(A') = \Derdc(A) + 
Res_{\overline{d} = \overline{c}} f_{A,A',\overline{d}}.
$$
\end{prop}
\Pr
To alleviate notations, we omit the evaluation at $a$ and the 
direction $\overline{d}$ in
the notations; to indicate summation along the arbitrary fixed 
direction $\overline{c}_{0}$, we just add the index $0$. Thus,
we respectively write:
\begin{eqnarray*}
F_{A} \text{~for~} S_{\overline{d}} \hat{F}_{A}(a) 
& \text{~and~} &
F_{A,0} \text{~for~} S_{\overline{c_{0}}} \hat{F}_{A}(a) \\
F_{A'} \text{~for~} S_{\overline{d}} \hat{F}_{A'}(a) 
& \text{~and~} &
F_{A',0} \text{~for~} S_{\overline{c_{0}}} \hat{F}_{A'}(a) \\
F_{A,A'} \text{~for~} S_{\overline{d}} \hat{F}_{A,A'}(a) 
& \text{~and~} &
F_{A,A',0} \text{~for~} S_{\overline{c_{0}}} \hat{F}_{A,A'}(a) \\
f_{A,A'} \text{~for~} f_{A,A',\overline{d}}
& \text{~and~} &
f_{A,A',0} \text{~for~} f_{A,A',\overline{c_{0}}}.
\end{eqnarray*}

>From the previous remark:
$$
F_{A'} = F_{A,A'} F_{A} \equiv 
(I_{n} + f_{A,A'}) F_{A} \pmod{\g^{\geq \delta}(\C)},
$$
so that:
$$
F_{A',0}^{-1} F_{A'}  \equiv 
F_{A,0}^{-1} F_{A,0} + f_{A,A'} - f_{A,A',0} \pmod{\g^{\geq \delta}(\C)}.
$$
The conclusion then comes by taking logarithms, applying the following
lemma and then taking residues.
\hfill $\Box$

\begin{lemma}
Let $M \in \G(\C)$ and $N \in \g^{\geq \delta}(\C)$. 
Then: 
$$
\log(M+N) \equiv (\log M) + N \pmod{\g^{\geq \delta}(\C)}.
$$
\end{lemma}
\Pr
Write $M = I_{n} + M'$. Then:
\begin{eqnarray*}
\log(M+N) 
& \equiv & \sum_{m \geq 1} \dfrac{(-1)^{m-1}}{m} (M' + N)^{m}  
\pmod{\g^{\geq \delta}(\C)}\\
& \equiv & \sum_{m \geq 1} \dfrac{(-1)^{m-1}}{m} {M'}^{m} + N 
\pmod{\g^{\geq \delta}(\C)}.
\end{eqnarray*}
\hfill $\Box$

\begin{cor}
Under the assumptions of the proposition, we have:
$$
\Derdc(A') - \Derdc(A) \in \g^{(\delta,\overline{c})}(\C).
$$
\end{cor}

We are going to prove that these are, in some sense, the only 
conditions on the $q$-alien derivations at a given level $\delta$.


\subsubsection{Interpolating categories}

There are two equivalent ways of defining $\EE^{(0)}_{p,1}$ from
$\EE^{(0)}_{1}$: the first is by restriction to a subclass of
objects, the pure ones; the second is by formalisation, \emph{i.e.}
extension of the base field $\Ka \rightarrow \Kf$. The former way
amounts to shrinking the Galois group $G^{(0)}_{1}$ to its \emph{quotient}
$G^{(0)}_{p,1}$. The latter way amounts to extending the class of
morphisms (indeed, there are no really new objects), and therefore
to shrinking the Galois group $G^{(0)}_{1}$ to its \emph{subgroup}
$G^{(0)}_{p,1}$. \\

The existence of a natural filtration on the Stokes group $\St$ suggests 
that it should be possible to interpolate between $\EE^{(0)}_{p,1}$ and 
$\EE^{(0)}_{1}$. We shall presently do so by extending the class of
morphisms; the interpretation by restriction to subobjects is a bit
more complicated. \\

We first define intermediate fields between $\Ka$ and $\Kf$,
for all levels $\delta \in \N$:
$$
\Kf^{(\delta)} = \left\{\sum f_{n} z^{n} \in \Kf \mid
\exists R > 0 \;:\; f_{n} = \text{O}(R^{n} q^{n^{2}/2 \delta})\right\}.
$$
Thus 
$\Kf^{(+\infty)} \underset{def}{=} \Ka \subset \Kf^{(\delta)} 
\subset \Kf^{(\delta - 1)} \subset \Kf^{(0)} \underset{def}{=} \Kf$
\footnote{
This $\Kf^{(\delta)}$ is the field of fraction of the algebra of
$q$-Gevrey series of level $\delta$, which was introduced
in \cite{Bezivin} and denoted $\C[[z]]_{q,s}$ with $s = 1/\delta$ 
(the $q$-Gevrey order) in \cite{RamisGrowth}.}.
The following is standard (\cite{RamisGrowth}, \cite{Zha02},\cite{JSSTO}):

\begin{lemma}
If $\nu - \mu = \delta \geq 1$, then, the following equation:
$$
(\sq F) (z^{\nu} B) - (z^{\mu} A) F = U, \quad
A \in GL_{r}(\C) \;,\; B \in GL_{s}(\C) \;,\; U \in \Ma_{r,s}(\Ka)
$$
has a unique solution $F \in \Ma_{r,s}(\Kf^{(\delta)})$. If moreover
$F \in \Ma_{r,s}(\Kf^{(\delta')})$ for some $\delta' > \delta$,
then $F \in \Ma_{r,s}(\Ka)$.
\end{lemma}

Write $\Kf^{> \delta} = \bigcup\limits_{\delta' > \delta} \Kf^{(\delta')}$.
Then we call $\Co^{\delta}$ the category with the same objects 
as $\EE^{(0)}_{1}$ (seen in matrix form) and with morphisms
satisfying the same conditions, but with $F \in GL_{n}(\Kf^{> \delta})$.
(Actually, since we deal only with integral slopes, we could as well
take $\Kf^{(\delta+1)}$ instead of $\Kf^{> \delta}$.) \\

It is then clear that the $\Co^{\delta}$ are tannakian categories,
and that embeddings are natural exact faithful $\otimes$-functors 
$\Co^{\delta} \rightarrow \Co^{\delta - 1}$. Moreover,
$\Co^{0} = \EE^{(0)}_{p,1}$, because equations with integral slopes
can be solved in $\Kf^{(1)}$ by the lemma, so that $\hat{F}_{A}$
is an isomorphism from $A_{0}$ to $A$ in $\Co^{0}$. In the opposite
direction, we have $\Co^{\infty} = \EE^{(0)}_{1}$. Actually,
if $A$ has slopes $\mu_{1} < \cdots < \mu_{k}$, then it is 
entirely determined by its image in $\Co^{\delta}$ for any
$\delta \geq \mu_{k} - \mu_{1}$. \\

>From the composite functor
$\Co^{\delta} \rightarrow \Co^{0} = \EE^{(0)}_{p,1}$, we draw that
objects in $\Co^{\delta}$ have a well defined Newton polygon, that
there is on $\Co^{\delta}$ a ``graded module'' functor,  and that
$\homa$ defines a fiber functor on $\Co^{\delta}$.

\paragraph{An alternative description of $\Co^{\delta}$.}

We also see that, if two objects $A$ and $B$ in $\EE^{(0)}_{1}$ have 
isomorphic images in $\Co^{\delta}$, then they have isomorphic images 
in $\Co^{0} = \EE^{(0)}_{p,1}$ and they can be written in standard form 
(\ref{equation:formestandard}) with the same block diagonal $A_{0}$. 
Of course, we may moreover assume $A_{0}$ to be in \emph{normalized} 
standard form and $A$, $B$ to be in \emph{polynomial} standard form
(section \ref{subsection:integralslopes}).

\begin{prop}
Let $A$ and $B$ in $\EE^{(0)}_{1}$ be in normalized polynomial standard 
form with the same block diagonal $A_{0}$. Then they have isomorphic images 
in $\Co^{\delta}$ if, and only if, there exists $F_{0} \in GL_{n}(\C)$
such that $B \equiv F_{0} A F_{0}^{-1} \pmod{\g^{> \delta}(\Ka)}$.
\end{prop}
\Pr
Here, of course, we have put
$\g^{> \delta} = \sum\limits_{\delta' > \delta} \g^{\delta'}$ (which is 
the same as $\g^{\geq \delta + 1}$ since we deal with integral slopes)
and the condition just means that $B$ and $F_{0} A F_{0}^{-1}$ have 
the same over-diagonals up to level $\delta$. \\

The diagonal part $F_{0}$ of any formal morphism $F$ from $A$ to $B$
is an automorphism of $A_{0}$, thus constant (because of normalisation).
Up to composing $F$ with $F_{0}^{-1}$ and replacing $A$ by 
$F_{0} A F_{0}^{-1}$, we may assume that $F_{0} = I_{n}$,
so that $F = \hat{F}_{A,B}$. The condition then means that $\hat{F}_{A,B}$
has its coefficients in $\Kf^{(\delta)}$. From the lemma, we draw, by 
induction on the level, that all over-diagonals up to level $\delta$
are analytic, therefore $0$ because of results in \cite{RSZ}.
\hfill $\Box$

\begin{cor}
One can define $\Co^{\delta}$ in the following alternative way:
\begin{enumerate}
\item{Objects of $\Co^{\delta}$ are matrices in $\EE^{(0)}_{1}$ modulo 
the equivalence relation $A \equiv B \pmod{\g^{> \delta}(\Ka)}$.}
\item{Morphisms from (the class of) $A$ to (the class of) $B$ are
matrices $F \in \Ma_{p,n}(\Ka)$ such that $(\sq F) A$ and $B F$ 
differ only in levels $> \delta$.}
\end{enumerate}
\end{cor}

\begin{cor}
The Galois group of $\Co^{\delta}$ is $\St(\delta) \rtimes G^{(0)}_{p,1}$
for some unipotent subgroup $\St(\delta)$ of $\St$. For $i \leq \delta$,
the $\Der^{(i)}_{\overline{c}}$ are well defined on $\Co^{\delta}$ and
belong to the Lie algebra $\st(\delta)$ of $\St(\delta)$.
\end{cor}


\subsubsection{A freeness theorem}

We now describe precisely the essential image of the functor $\F$,
that is, given $A_{0}$ in $\EE^{(0)}_{p,1}$, the exact conditions 
on $(\Derdc(A))_{\delta \geq 1,\overline{c} \in \Eq}$ that allow
the reconstruction of $A$. The reconstruction will be done inductively,
using $q$-alien derivations of levels up to $\delta$ to reconstruct
the over-diagonals of $A$ up to level $\delta$, that is (after the
previous paragraph) an object in $\Co^{\delta}$. \\

In the same spirit as the definition of isoformal analytic classes 
in \cite{RSZ}, we consider classes of objects $A$ in $\Co^{\delta}$
above an object $B$ of $\Co^{\delta-1}$ under the equivalence
induced by gauge transform $F \equiv I_{n} \pmod{\g^{\leq \delta}}$.
Using polynomial standard normal form and the results of \emph{loc. cit.},
we see that these classes make up a vector space of dimension:
$$
\text{irr}^{\delta}(A_{0}) = 
\sum_{\mu_{j} - \mu_{i} = \delta} r_{i} r_{j} (\mu_{j} - \mu_{i}) =
\delta \sum_{\mu_{j} - \mu_{i} = \delta} r_{i} r_{j}.
$$
Moreover, to see if two objects $A,A'$are in the same class, one 
computes $\hat{F}_{A,A'} \in \G^{\delta}(\Kf)$; if its over-diagonal
at level $\delta$ has null $q$-Borel invariants, then we have the same
class. 

\begin{thm}[Freeness theorem]
Let $B$ be an object of $\Co^{\delta-1}$. Then, there is an affine
space $V_{\overline{c}}(B)$ of direction $\g^{(\delta,\overline{c})}(\C)$ 
such that: \\
(i) The $\Derdc(A)$ for $A$ an object of $\Co^{\delta}$ above $B$
belong to $V_{\overline{c}}(B)$. \\
(ii) The mapping which sends an object $A$ of $\Co^{\delta}$ above $B$
to the family of all $\Derdc(A)$ induces a one-to-one correspondance 
between classes of such objects (as defined above) and 
$\prod V_{\overline{c}}(B)$.
\end{thm}
\Pr
(i) It follows from paragraph \ref{subsubsection:arbitrarylevel} that 
all $\Derdc(A)$ (where $A$ is fixed and $\overline{c}$ varies) belong 
to a unique affine space of direction $\g^{(\delta,\overline{c})}(\C)$.
Call it $V_{\overline{c}}(B)$. It is easily seen that the product space 
$\prod V_{\overline{c}}(B)$ has dimension $\text{irr}^{\delta}(A_{0})$. \\
(ii) The map from the set of polynomial representatives of a class, 
as described above, onto the above affine space, is affine.
Up to the choice of an arbitrary basepoint, it is equivalent, after
the results of \cite{RS3} (section 3.2), to the parametrisation of
the isoformal class by $q$-Borel transform, which is one-to-one after
\cite{RSZ}.
\hfill $\Box$

\begin{cor}
The following algorithm allows one to reconstruct $A$ in $\EE^{(0)}_{1}$ 
from $A_{0}$ and the $\Derdc(A)$:
\begin{enumerate}
\item{Reconstruct the first over-diagonal using the $q$-derivations of
lowest level (this is the linear situation and it rests on \cite{RS1}.}
\item{Having reconstructed the over-diagonals up to level $\delta-1$
(using $q$-alien derivations up to level $\delta - 1$), call $A'$ the
matrix with these over-diagonals and $0$ above; then compute the $\Derdc(A')$.}
\item{Use the relation 
$\Derdc(A') - \Derdc(A) \in \g^{(\delta,\overline{c})}(\C)$
to find the level $\delta$ over-diagonal of $\hat{F}_{A,A'}$,
then the level $\delta$ over-diagonal of $A$.}
\end{enumerate}
\end{cor}

In \cite{RS3}, we shall give a representation-theoretic formulation
of the theorem, and a description of the nonlinear part of $\Derdc(A)$
(that part which depends on the lower level $q$-derivations) in terms
of convolution.



\section{Conclusion}
\label{section:conclusion}

Write $H = \C^{*} \times Hom_{gr}(\C^{*}/q^{\Z},\C^{*})$ and
$\nu = (t,\gamma) \in H$. We saw how the group $H$ acts upon 
the diagonal of $A_{0}$: for $1 \leq i \leq k$ and $\alpha \in \Sp(A_{i}$,
positions corresponding to slope $\mu_{i}$ and the eigenspace of $A_{i,s}$ 
for $\alpha$ are multiplied by $t^{\mu_{i}} \gamma(\alpha)$. Now let
$i < j$ be indices of slopes $\mu_{i} < \mu_{j}$ and $\alpha \in \Sp(A_{i}$,
$\beta \in \Sp(A_{j})$ be corresponding exponents. The adjoint action of $\nu$
on the $(\mu_{i},\overline{\alpha},\mu_{j},\overline{\beta})$ block is the
multiplication by:
$$
t^{\mu_{i} - \mu_{j}} \dfrac{\gamma(\alpha)}{\gamma(\beta)} =
\left(t^{-1} \gamma(\overline{c})\right)^{\delta}
$$
for each ``resonant'' $\overline{c}$, \emph{i.e.} $\overline{c} \in \Eq$
such that $c^{\delta} \equiv \dfrac{\alpha}{\beta} \pmod{q^{\Z}}$. 
For any Galois derivation $D \in \st$, we now put:
$$
\Phi_{\nu}^{(\delta,\overline{c})}(D) = 
\nu(D) - \left(t^{-1} \gamma(\overline{c})\right)^{\delta} D \in \st,
$$
where $\nu(D)$ comes from the adjoint action of $H$ on $\st$.
>From the remark on page \pageref{remark:espacefixe} and from 
paragraph \ref{subsubsection:arbitrarylevel}, one draws that,
for two objects $A,A'$ of $\Co^{\delta}$ above the same object 
$B$ of $\Co^{\delta-1}$, 
$\Phi_{\nu}^{(\delta,\overline{c})}(\Derdc)(A)
\Phi_{\nu}^{(\delta,\overline{c})}(\Derdc)(A') = 0$. In other
words, $\Phi_{\nu}^{(\delta,\overline{c})}(\Derdc)(A)$ depends 
only on the lower levels $\delta' < \delta$ of $A$. Moreover
it is trivial on the first level. Actually, with methods similar
to those used here, one can prove that 
$\Phi_{\nu}^{(\delta,\overline{c})}(\Derdc)(A)$
is in the Lie algebra generated by the $q$-alien derivations
at lower levels. So it is natural to conjecture that 
$\Phi_{\nu}^{(\delta,\overline{c})}(\Derdc)$ belong to the free 
Lie algebra generated by the $\Der^{(\delta')}_{\overline{d}}$
($\delta' < \delta$, $\overline{d} \in \Eq$), and even that there
is a universal explicit formula. This would allow us to define
a semi-direct product by a free Lie algebra, and to definitely
``free'' the $q$-alien derivations. \\

All the problems comes from the fact that points come from two
distinct origins: elements of the dual of $H$ on the one hand, 
packs of points of $\Eq$ on the other hand, and from the interplay
of the corresponding games of localisation. Comparing with the
differential case, where one localizes geometrically on the circle 
of directions $S^{1}$, \emph{then} one takes a Log, here, we take 
a Log, \emph{then} we localise on $Eq$; whence an embroilment
with plenty of Campbell-Haussdorff formulas between the two approaches
\footnote{Actually, we think that, in the end, we'll have a simler
description with a denumerable family of $q$-alien derivations,
freed by the mere action of the theta torus.}. \\

We shall also give in \cite{RS3} various applications, to the abelianisation
of the tannakian $\pi_{1}$ and to the inverse problem for the local Galois
group. For the latter problem, we shall state a list of necessary conditions;
we don't know for the time being if they are sufficient. \\

Last, we built our alien $q$-derivations by tannakian methods. One can ask
what happens for \emph{solutions}. There, one meets the usual difficulty 
about constants, since one wishes operators defined over $\C$ and acting
upon solutions (while constants are here elliptic functions). That problem
maybe has no solution; however, one could perhaps, in analogy with the
differential case, ``unpoint'' the $q$-alien derivations and build operators
acting upon adequate spaces of formal power series. This seems related to
a ``$q$-convolution'' mechanism presently studied by Changgui Zhang.

\vskip 10pt

\noindent{\large\bf Acknowledgements}

\vskip 5pt

The work of the first author has been partially supported 
by the NEST EU Math. Project GIFT, FP6-005006-2.



\end{document}